\documentclass{elsarticle}
\usepackage{caption}
\usepackage{epsfig,multirow,bbm}
\usepackage{cancel}
\usepackage{afterpage}
\usepackage{gensymb}
\usepackage{amsmath}
\usepackage{amssymb}
\usepackage{dsfont}
\usepackage{url,color}
\usepackage{booktabs}
\usepackage{endnotes}
\usepackage{subfigure}

\graphicspath{{./}}
\DeclareMathAlphabet{\mathpzc}{OT1}{pzc}{m}{it}

\newcommand{\dt}{\Delta t}
\newcommand{\dx}{\Delta x}

\newcommand{\ones}{\mathbbm{1}}

\newcommand\ST{\rule[-0.55em]{0pt}{1.8em}}

\newcommand{\mratio}{m}

\usepackage{theorem}

\theorembodyfont{\upshape}

\begin{document} 
\begin{frontmatter}

\title{Implicit Extensions of an Explicit Multirate Runge--Kutta Scheme}

\author{Emil M. Constantinescu}
\ead{emconsta@mcs.anl.gov}

\address{Mathematics and Computer Science Division, Argonne National
  Laboratory, \\9700 S. Cass Avenue, Argonne, IL 60439, USA, Tel. +1 630
  252 0926}

\begin{abstract}
We propose a new method that extends conservative explicit multirate
methods to implicit explicit-multirate methods. We develop extensions
of order one and two with different stability properties on the
implicit side. The method  is suitable for time-stepping
adaptive mesh refinement PDE discretizations with different degrees
of stiffness. A numerical example with an advection-diffusion problem
illustrates the new method's properties.
\end{abstract}


\end{frontmatter}

%
\section{Introduction}
%

Multirate time-stepping methods have become popular in solving
computational fluid dynamics problems with adaptive meshes
\cite{Mittal_2021,Roberts_2021,Grote_2015}. Explicit
multirate methods are particularly efficient at integrating partial
differential equations (PDE)s with
adaptive mesh refinement (AMR), where local time stepping overcomes the
global Courant--Friedrichs--Lewy (CFL) limitation \cite{Seny_2013}. If these problems
have stiff components, however, using purely
explicit methods is inefficient or impractical. This situation forces one to use implicit-explicit (IMEX) methods; however, existing IMEX conservative methods treat the explicit component
with the same global time step, which limits the time step due to the fastest
component---typically, associated with the finest grid
points or elements. This restriction leads to ineffective time stepping
especially  when one uses AMR for stiff problems. While general
frameworks have been proposed (e.g., \cite{Sandu_2019}), little work has been done on methods that are suitable for such conditions
and also preserve linear invariants.

In this work we propose an extension to explicit multirate methods
with a computationally efficient implicit component, which allows
multirate treatment of the explicit component and implicitness for
stiff components. The stiff components are treated with a single rate.   
We consider the following initial value
problem: 
$\dot{y}(t)=F(y(t))=f(y(t)) + g(y(t)) \,,~ t_0\le t \le
t_F\,,~ y(t_0)=y_0\,. $
Here $y\in\mathbb{R}^N$, and $F$ typically represents the spatial
discretization and can be additively partitioned into
$f,g:\mathbb{R}^N\rightarrow\mathbb{R}^N$, which are Lipschitz continuous
functions. We assume that $f$ is a nonstiff component that can be
efficiently integrated with an explicit integrator and $g$ is a
component that can be stiff and requires an implicit integrator.
We further assume that the solution $y$ can be partitioned
in non-overlapping subdomains, where the dynamic behavior with respect to $f$ is different. The explicit stability
requirements for each subdomain are separated in different classes
called levels---a reference to AMR. Therefore, we consider a component
partitioning 
$y^\top=[y^F, y^S]$ and the same for the explicit components $f^\top=[f^F,
  f^S]$ so that the partitioned system to be solved takes the
following form,
%
\begin{align}
\label{eq:ODE:partitioned}
\begin{bmatrix}
  \dot{y}^F \\
  \dot{y}^S
\end{bmatrix}
=
\begin{bmatrix}
  f^F\left(y^F,y^S\right)  \\
  f^S\left(y^F,y^S\right) 
\end{bmatrix}
+
  g\left(y\right)
\,,
\end{align}
%
where the superscript $F$ refers to the fast component  and superscript $S$ to the slow
component. For brevity we consider only two partitions
here; but, in general, multiple partitions can be accommodated, which are
relevant for practical multiple-level adaptive mesh refinement.  

A class of conservative multirate explicit partitioned Runge--Kutta (MPRK) methods was
introduced in \cite{Constantinescu_A2007e}. These methods are
constructed by using a base method and can attain second-order overall
accuracy if the base method is at least second order.  In this work we
extend a particular MPRK method by adding an implicit stage to handle stiff
additive partitions along with the original explicit multirate scheme,
while maintaining conservation at the temporal discrete level. In the
PDE case, so long as the spatial discretization is conservative, the
time-marched solution with the proposed method is also conserved.

The rest of the paper is organized as follows. Section \ref{sec:MPRK}
introduces MPRK and the notation. We provide the new method design
and specific examples in Sec.~\ref{sec:MPIMEXRK}. We illustrate the properties of the
newly introduced methods in Sec. \ref{sec:numerical} through numerical
examples. We provide concluding
remarks in Sec. \ref{sec:conclusions}.

%
\section{MPRK Methods\label{sec:MPRK}}
%

The multirate explicit partitioned Runge--Kutta scheme as 
introduced in \cite{Constantinescu_A2007e} is constructed by using a \emph{base}
method that is repeated $\mratio$ times with a fractional time step. The base method is defined as 
a classical RK method by coefficients $A=[a_{i,j}]$, $b=[b_i]$, and
$c=[c_i]=A \ones_s$, where $\ones_s$ is a vector of ones of length
$s$ and $i,j=1,\dots,s$ and represented in a tableau, 
%
$  \begin{array}{c|c}
\ST    c & A\\
\hline
\ST     & b^\top
\end{array}$. %
In order to preserve numerical stability and conservation properties,
the \emph{fast} or the subcycled ($\mratio$ times) component that 
corresponds to the fast partition is a repeated application of the
base method with a time step of $\frac{\dt}{\mratio}$; we refer to
$\mratio$ as the multirate ratio. 
The region around the boundary between fast and slow regions is called
a \emph{buffer}, with a size that depends on the spatial discretization (stencil and rate), and is
typically small relative to the fast and slow regions. 
In the buffer region one simply
repeats the base method $\mratio$ times with $\dt$,  leading to the
\emph{slow} method so that both the slow and the fast methods 
have the same number of stages. Then, inside the slow
region, the slow method simply reverts to the base method. Efficiency
gains result from applying a small time step on the fast regions and large
steps on the slow regions. By construction, the fast method has a factor
of $\mratio$ more stages  than the
slow method on the slow method has and a fraction step size, $\dt/\mratio$.  
For multiple refinement levels, one can use
telescoping nesting, which is achieved by replacing the base method
with the fast method and repeating the
procedure above.

A fast-slow method can be represented as a partitioned RK:
\begin{subequations}
  \label{eq:ODE:partitioned:MPRK}
  \begin{align}
    \label{eq:ODE:partitioned:MPRK:fast}
Y^F_i=& y^F_n + \dt \sum\nolimits_{j=1}^{i-1} a^F_{ij}
f^F\left(Y^F_j,Y^S_j\right) \,,~i=1,\dots,s\\
\label{eq:ODE:partitioned:MPRK:slow}
Y^S_i=& y^S_n + \dt \sum\nolimits_{j=1}^{i-1} a^S_{ij}
f^S\left(Y^F_j,Y^S_j\right) \,,~i=1,\dots,s \\
\label{eq:ODE:partitioned:MPRK:comp}
y_{n+1} =& y_n+ \dt \sum\nolimits_{i=1}^{s} b_i f(Y_i)  \,,
\end{align}
\end{subequations}
where the stages are partitioned as the solution in two fast and slow components,
$Y_i^\top = [Y^F_i,Y^S_i]$ and
$f(Y_i)^\top=[f^F\left(Y^F_i,Y^S_i\right),f^S\left(Y^F_i,Y^S_i\right)]$.
We also set $b^F=b^S=b$, which is required
for conservation. This results in the short form of
\eqref{eq:ODE:partitioned:MPRK:comp}, where each component is
computed by $y_{n+1}^{\{F,S\}} = y_n^{\{F,S\}}+ \dt
\sum\nolimits_{i=1}^{s} b_i^{\{F,S\}} f^{\{F,S\}}(Y^F_i,Y^S_i)$. The
step size is controlled by the $a^{\{F,S\}}$ coefficients.

An example with $\mratio=2$ is
\begin{align}
  \label{eq:Butcher:MPRK2}
  \begin{array}{c|cc}
\ST    0 & 0\\
\ST    1 & 1 &0\\
\hline
\ST     & \frac{1}{2} & \frac{1}{2}\\
  \end{array}
  \qquad
\begin{array}{c|cccc}
\ST    0 & 0\\
\ST    1 & 1 &0\\
\ST    0 & 0 &0 & 0\\
\ST    1 & 0 &0 & 1 &0\\
\hline
\ST     & \frac{1}{4} & \frac{1}{4} & \frac{1}{4} &\frac{1}{4}  \\
  \end{array}
  \qquad
    \begin{array}{c|cccc}
\ST    0 & 0\\
\ST    \frac{1}{2} & \frac{1}{2} &0\\
\ST    \frac{1}{2} & \frac{1}{4} &\frac{1}{4} &0\\
\ST    1 & \frac{1}{4} &\frac{1}{4} & \frac{1}{2}  &0\\
\hline
\ST     & \frac{1}{4} & \frac{1}{4} & \frac{1}{4} &\frac{1}{4}  \\
\end{array}\,,
\end{align}
where the tableaux represent the base ($a_{ij},b_{i},c_{i}$), the slow method
($a^S_{ij},b^S_{i},c^S_{i}$), and the fast method
($a^F_{ij},b^F_{i},c^F_{i}$). In the fast region, the 4-stage fast method is used; and in
the buffer and slow regions, the slow method is used. In the
slow region, however, the slow method is equivalent to the 2-stage base
method, and hence the computational savings is relative to having a global
time step. The fast method is the result of applying the base method
twice with a step size of $\dt/2$; therefore the fast method takes
takes twice as many half-steps as the slow method takes.
This particular method is the focus of our study here.
The reader is referred to \cite{Constantinescu_A2007e} and
a  recent algorithmic presentation in \cite{Kang_2021b} for a more
in-depth description of these methods.

%
\section{Extended MPRK-Implicit Methods\label{sec:MPIMEXRK}}
%

We now introduce a three-way partitioning of
\eqref{eq:ODE:partitioned:MPRK} into explicit fast, explicit slow, and
implicit stiff. We restrict the presentation to two levels; however,
the explicit component of the method
can be extended to arbitrary nesting levels in the same fashion as discussed
above, while maintaining a fixed ratio. The resulting scheme is now
described by three sets of 
coefficients: superscripts $F$, $S$, and $\tilde{}$ define the fast,
slow, and  implicit coefficients, respectively. We have two goals:
($i$) minimize the number of implicit stages and ($ii$) accommodate
multiple levels of refinement, that is, telescoping nesting. To
accomplish both and to avoid the loss of conservation, we will modify
only the last stage in \eqref{eq:ODE:partitioned:MPRK}, which will be
common across all partitions. 
Therefore, the last stage and the step completion in
\eqref{eq:ODE:partitioned:MPRK} become   
\begin{subequations}
  \label{eq:ODE:partitioned:RK}
\begin{align}
Y_s=& y_n + \dt \sum_{j=1}^{s-1} \begin{bmatrix}  a^F_{sj}
  f^F\left(Y_j\right) \\ \ST a^S_{sj} f^S\left(Y_j\right) \end{bmatrix} +
\dt \sum_{j=1}^{s} \tilde{a}_{sj}
g\left(Y_j\right) \\
y_{n+1} =& y_n+ \dt \sum_{i=1}^{s} b_i \left(f(Y_i) + g(Y_i) \right)  \,.
\end{align}
\end{subequations}
The implicit component will have only one implicit
stage regardless of the number of explicit partitions and of the ratio among
them. Therefore, the methods introduced here are neither multirate implicit  \cite{Carciopolo_2019} nor multirate IMEX \cite{Constantinescu_2013}. The
Butcher tableau representation of the implicit component indexed by
$\tilde{a}$ as a matrix will therefore be all zeros except in the last stage and can
augment the tableaux in \eqref{eq:Butcher:MPRK2}. With the $b$ vector
fixed, the choice of $\tilde{a}_{sj}$, $j=1,\dots,s$ is driven by stability and
lastly by accuracy considerations.

We analyze stability by considering the following problem, $\dot{y}=
 \lambda^F y + \lambda^S y + \lambda^\textnormal{stiff}y$, with
 $\lambda^{\{F,S,\textnormal{stiff}\}} \in \mathbb{C}$ 
resulting from linear operators that can be diagonalized simultaneously, where
each component on the right-hand side is time-stepped with the respective fast,
slow, and implicit integrators. The stability function is defined by
 $R(z)$, $z=\lambda \dt$ that satisfies $y_{[n+1]} = R(z)
y_{n}$. Individually, the stability of the fast method is $\mratio$
times larger than that of the slow method. More complex stability 
analyses are possible, in particular if interpolation is
needed \cite{Bonaventura_2020}. In general, simplifying assumptions have
to be made in all cases; however, they tend to hold well in practice.

%
\subsection{A-Stable Second-Order Extensions}
\label{sec:A-Stable}
%
A multirate-explicit A-stable-implicit method 
that is second-order explicit and second-order
implicit A-stable is obtained by setting $\tilde{a}_{sj}=\frac{1}{2}$,
$j=1,\dots,s$ in \eqref{eq:ODE:partitioned:RK}.
The resulting
stability function of the implicit method is $R(z) = \frac{2+z}{2-z}$,
which implies that the implicit part is A-stable with
$R(\infty)=-1$. By using these coefficients we guarantee that the
implicit-explicit  stability function damps large eigenvalues;
see 
\cite[Sec 4.5]{Giraldo_2012P}. The abscissa of the implicit method becomes  
$\tilde{c}=[0,0,\dots,0,\mratio]^\top$ on the fast and slow parts and
$\tilde{c}=[0,1]^\top$ for the base method. In both cases $\tilde{c}
b=\frac{1}{2}$, which gives second order of the implicit part and of
the overall method. Computationally the implicit part has only one implicit
stage, with the rest of the stage coefficients being zero. The
resulting method can be expressed in 
the same fashion as MPRK by augmenting the following tableaux to
\eqref{eq:Butcher:MPRK2}: 
\begin{align}
  \label{eq:Butcher:MPRK2:IM:AStable:2}
   \begin{array}{c|cc}
\ST    0 & 0\\
\ST    1 & \frac{1}{2} & \frac{1}{2}\\
\hline
\ST     & \frac{1}{2} & \frac{1}{2}\\
  \end{array}
  \qquad
    \begin{array}{c|cccc}
\ST    0 & 0\\
\ST    0 & 0 & 0\\
\ST    0 & 0 & 0 &0\\
\ST    2 & \frac{1}{2} &  \frac{1}{2} & \frac{1}{2} & \frac{1}{2}\\
\hline
\ST     & \frac{1}{4} & \frac{1}{4} & \frac{1}{4} &\frac{1}{4}  \\
\end{array}\,,
\end{align}
where the first tableau corresponds to the base method and the second
tableau to the slow and fast methods in \eqref{eq:Butcher:MPRK2}. In
particular such a method corresponds to the following algorithm:
%
\begin{align*}
  Y^F_1=& y^F_n \,, \qquad Y^S_1= y^S_n\\
  Y^F_2=& y^F_n + \frac{\dt}{2}  f^F\left(Y^F_1,Y^S_1\right)\,, \qquad
  Y^S_2= y^S_n + \dt f^S\left(Y^F_1,Y^S_1\right) \,, \\
  Y^F_3=& y^F_n + \frac{\dt}{4}  f^F\left(Y^F_1,Y^S_1\right) + \frac{\dt}{4}  f^F\left(Y^F_2,Y^S_2\right)\,, \qquad
  Y^S_3= y^S_n \,, \\
{\color{red}Y_4}=& y_n + \dt \begin{bmatrix}   \frac{1}{4}  f^F\left(Y^F_1,Y^S_1\right) +
  \frac{1}{4}  f^F\left(Y^F_2,Y^S_2\right)+ \frac{1}{2}
  f^F\left(Y^F_3,Y^S_3\right)\\ \ST
  f^S\left(Y^F_3,Y^S_3\right)  \end{bmatrix} +\\
& \quad \frac{\dt}{2} \left( g(Y_1)+g(Y_2)+g(Y_3)+g({\color{red}Y_4})\right)\,, \\  
  y_{n+1} =& y_n+ \frac{\dt}{4}\sum_{i=1}^{4}\left(f(Y_i)+g(Y_i)\right) \,,
\end{align*}
%
where the terms in red illustrate the sole implicit components.

%
\subsection{L-Stable First-Order Extensions}
\label{sec:L-Stable}
%

Partitioned methods of type
\eqref{eq:ODE:partitioned:MPRK}, \eqref{eq:ODE:partitioned:RK}  with
one implicit stage cannot be both second order and L-stable at the
same time because of the constraints imposed on the order conditions
and the stability function. Thus, L-stable implicit 
methods coupled to the multirate method can be at
most second order accurate on the explicit part and first order on the implicit
part. One such method is obtained by choosing $\tilde{a}_{sj}=1$,
$j=1,\dots,s$ in \eqref{eq:ODE:partitioned:RK}. The resulting
stability function of the implicit method is $R(z) = \frac{1}{1-z}$. The
abscissa of the implicit method corresponds to
$\tilde{c}=[0,0,\dots,0,\mratio\,s]^\top$ on the fast and slow parts and
$\tilde{c}=[0,s]^\top$ for the base method and  $\tilde{c}
b=\frac{1}{2}=1$. Therefore, this multirate-IMEX method is
second-order explicit and first-order 
implicit L-stable.  As in the case above, the resulting method can be
expressed in 
the same fashion as MPRK by augmenting the following tableaux to \eqref{eq:Butcher:MPRK2}:
\begin{align}
  \label{eq:Butcher:MPRK2:IM:LStable:1}
  \begin{array}{c|cc}
\ST    0 & 0\\
\ST    2 & 1 & 1\\
\hline
\ST     & \frac{1}{2} & \frac{1}{2}\\
  \end{array}
  \qquad
  \begin{array}{c|cccc}
\ST    0 & 0\\
\ST    0 & 0 & 0\\
\ST    0 & 0 & 0 &0\\
\ST    4 & 1 & 1 &1 &1\\
\hline
\ST     & \frac{1}{4} & \frac{1}{4} & \frac{1}{4} &\frac{1}{4}  \\
\end{array}\,,
\end{align}
where the first tableau corresponds to the base method and the second
tableau to the slow and fast methods in \eqref{eq:Butcher:MPRK2}.

\newpage
%
\section{Numerical Examples\label{sec:numerical}}
%
We illustrate the properties of the methods introduced above by using
a one-dimensional advection-diffusion problem,

\noindent
\begin{minipage}{0.54\textwidth}
\begin{align}
  \label{eq:PDE:advection:diffusion}
\frac{\partial{u(t,x)}}{\partial{t}} +
\frac{\partial}{\partial{x}} \left(\omega(x)u(t,x) \right) = \delta\, \frac{
  \partial^2{u(t,x)}}{\partial{x}^2}\,,  
\end{align}
\null
\par\xdef\tpd{\the\prevdepth}
\end{minipage}
\begin{minipage}{0.44\textwidth}
\centering
\includegraphics[width=0.80\textwidth]{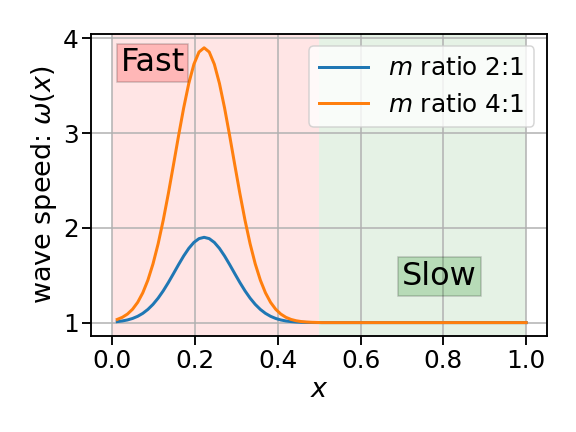}
\vspace{-6mm}
\captionof{figure}{Flux function ($\mratio=2$ \&
  $\mratio=4$)\label{fig:flux}}
\end{minipage}

\noindent{}which is discretized in space in $M$ uniform intervals. This leads to a
semi-discrete system $y(t) = \{u_k(t)\}_{k=1,\dots,M}$.
%
%
We discretize the advective term by using a conservative third-order
upwind-biased finite difference method and the diffusive term by
second-order finite differences, and we use periodic boundary conditions. The
numerical stability restriction in the advective term is adjusted by
using different values of $\omega$ in space to mimic AMR
(Fig. \ref{fig:flux}), where we consider two cases suitable for $\mratio=2$ and
$\mratio=4$. The stiffness is controlled by changing the value of
$\delta$. The time step is kept the same for all simulations.
Advection is treated explicitly with different rates and diffusion
implicitly with a single-rate integrator.

\begin{figure}[!ht]
\begin{center}
Mildly stiff $\rightarrow$ $\delta=.05$; mass loss: \ref{fig:mrk2}: 2.4e+50;
\ref{fig:srk2}: 0.0 ; \ref{fig:mrk2:imex2}: 1.1e-16
\\
\subfigure[\vspace{-0.4mm}explicit multirate \eqref{eq:Butcher:MPRK2}]{\includegraphics[width=0.32\textwidth]{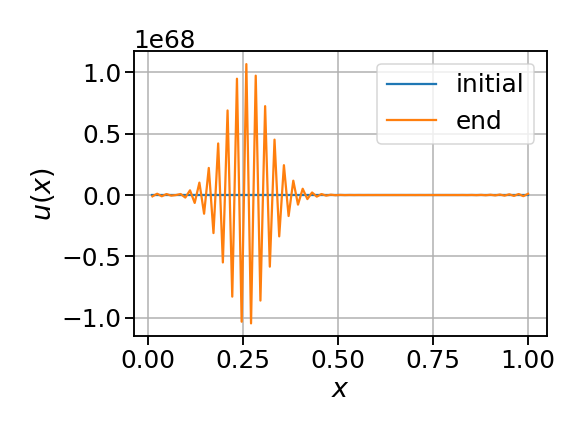}\label{fig:mrk2}}
\subfigure[single-rate IMEX]{\includegraphics[width=0.32\textwidth]{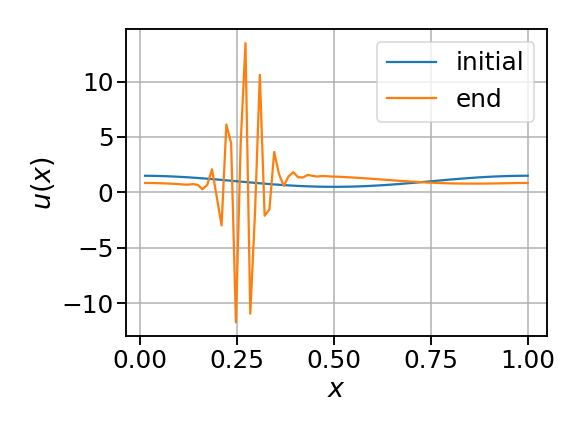}\label{fig:srk2}}
\subfigure[multirate w/ implicit \eqref{eq:Butcher:MPRK2:IM:AStable:2}]{\includegraphics[width=0.32\textwidth]{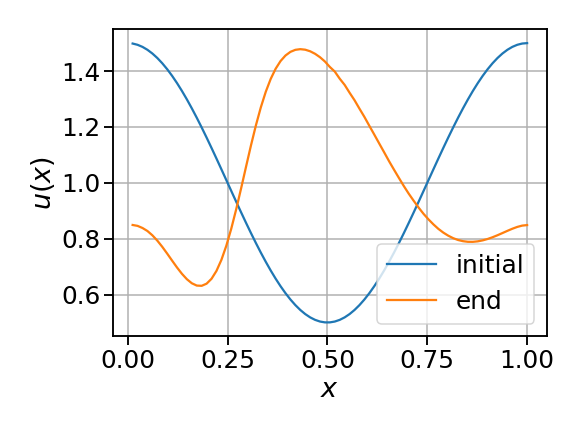}\label{fig:mrk2:imex2}}
\\
Stiff $\rightarrow$ $\delta=100$; mass loss:
(\ref{fig:mrk2:imex2:high:delta}: 4e-13;
\ref{fig:mrk2:imex1:high:delta}: 6e-13; \ref{fig:mrk2:imex2:m4}: 7.8-16
\\
\subfigure[A-stable multirate \eqref{eq:Butcher:MPRK2:IM:AStable:2}]{\includegraphics[width=0.32\textwidth]{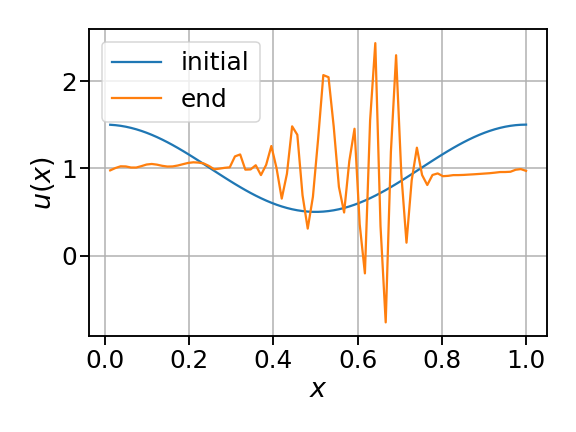}\label{fig:mrk2:imex2:high:delta}}
\subfigure[L-stable multirate \eqref{eq:Butcher:MPRK2:IM:LStable:1}]{\includegraphics[width=0.32\textwidth]{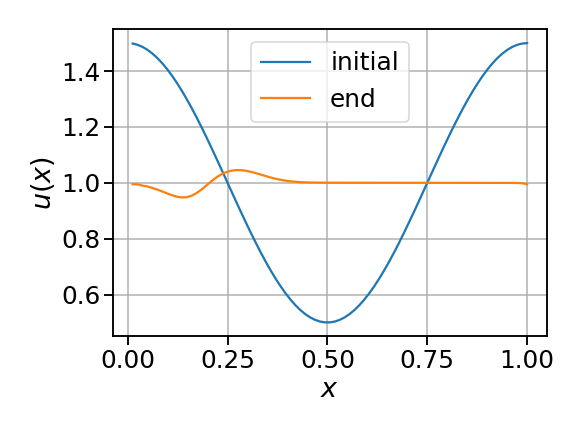}\label{fig:mrk2:imex1:high:delta}}
\subfigure[\ref{fig:mrk2:imex2} with $\mratio=4$]{\includegraphics[width=0.32\textwidth]{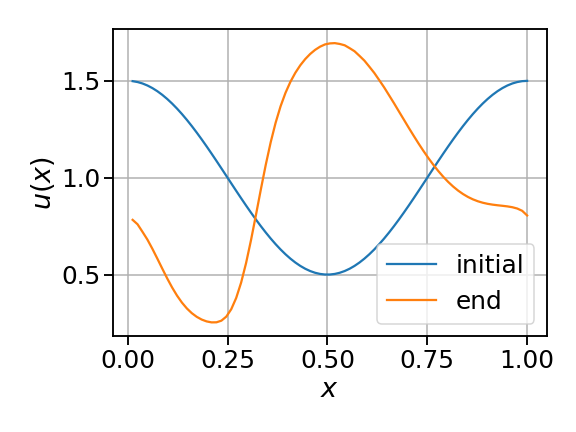}\label{fig:mrk2:imex2:m4}}
\end{center}
\caption{Results integrating \eqref{eq:PDE:advection:diffusion} by
  using a finite volume method in space and various time integrators
  with $\dx=0.012$ ($M=81$ grid points), fixed $\dt=0.0125$ (CFL=$\omega(x)$), and final
  time $t_F=0.3$. The resulting advective CFL is 1.92 on the fast region and 1.01 on
  the slow region. The discrete mass loss is defined as $\dx \left\|
  \sum_{k=1}^{M} u_k(t_0) - \sum_{k=1}^{M} u_k(t_F) \right\| $.\label{fig:results}} 
\end{figure}

In the case of mildly stiff problems, $\delta=0.05$, the original
MPRK \cite{Constantinescu_A2007e} solution with $\mratio=2$ in
Fig. \ref{fig:mrk2} is unstable because of problem stiffness, and the
single-rate IMEX solution is 
unstable because the explicit part violates local CFL conditions
(Fig. \ref{fig:srk2}). The 
multirate ($\mratio=2$) explicit-implicit A-stable solution combines
the stability features of 
both methods and is stable and efficient because of
local time stepping on the explicit part and the implicit component
(Fig. \ref{fig:mrk2:imex2}). In Fig. \ref{fig:mrk2:imex2:m4} we show
the solution  in the same setting as shown in Fig. \ref{fig:srk2}
except that the wave speed takes a 4:1 ratio between fast and slow
(see Fig. \ref{fig:flux}) and we use an $\mratio=4$, which means that
the fast integrator takes four time steps with $\dt/4$ for each slow
integrator step in the slow region. 

For stiff problems,
$\delta=100$, the A-stable method (\S\ref{sec:A-Stable}) is
not stable (see Fig. \ref{fig:mrk2:imex2:high:delta}); however, the
L-stable method (\S\ref{sec:L-Stable}) becomes stable (see
Fig. \ref{fig:mrk2:imex1:high:delta}) at the cost of reducing the
implicit accuracy order to one. 
%

%
\section{Conclusions\label{sec:conclusions}}
%
%
We propose a new conservative method that extends an explicit multirate
method to A- and L-stable implicit multirate-explicit methods. These
extensions are first and second order and can accommodate telescopic
multirate nesting. Therefore, these methods are suitable for AMR PDE
discretizations with stiff terms that are being treated with a single-rate method. Such time-stepping schemes remove
the stiffness stability limitation of existing conservative explicit multirate methods 
applied to adaptive mesh refinement PDE discretizations,  thus
allowing them to be applied to problems with  stiff
components as well. An advection diffusion problem is used to demonstrate the
stability properties of the new methods.

\section*{Acknowledgments}
\noindent{}{\small{}This material is based upon  work  supported by the U.S. Department of
Energy, Office of Science, Office of Advanced Scientific Computing
Research (ASCR) program and through the Fusion Theory Program of the
Office of Fusion Energy Sciences and the SciDAC partnership on Tokamak
Disruption Simulation between the Office of Fusion Energy Sciences and
the Office of ASCR, under Contract DE-AC02-06CH11357. I would also like to
thank the reviewers for their valuable input.}

\end{document}